\documentclass{amsart}
\usepackage{txfonts}
\usepackage{graphics}
\theoremstyle{plain}
\newtheorem{theorem}{Theorem}
\newtheorem{corollary}{Corollary}
\begin{document}
\title{On the decycling of powers and products of cycles}
\author{Adrian Riskin}
\address{Department of Mathematics\\
Mary Baldwin College\\
Staunton, Virginia 24401 USA}
\email{ariskin@mbc.edu}
\thanks{I would like to thank Sarah Lynch and her colleagues for providing a congenial
work environment.  Also I would like to thank my colleague Martha Walker for her interest
and encouragement.}
\begin{abstract}
We calculate exact values of the decycling numbers of $C_{m} \times C_{n}$ for $m=3,4$, 
of $C_{n}^{2}$, and of $C_{n}^{3}$.
\end{abstract}
\maketitle

\section{Introduction and definitions}

If $G$ is a simple graph, $S \subseteq V(G)$, and $G-S$ is acyclic, then $S$ is a 
\textit{decycling set} of $G$.  The size of a smallest decycling set is the 
decycling number $\nabla(G)$ of $G$.  It turns out that, in contrast to the corresponding
problem for edges, finding the decycling number can be a quite difficult problem, even for
some very simple families of graphs.  The problem is NP-complete in general [4].  Much 
recent work in the area has been focused on calculating the decycling numbers of various
families of graphs, e.g. hypercubes, grids, and so on.  For a good introduction to and 
bibliography of recent work of this nature see [1,3].  Most of these results, rather than
yielding exact values, are instead in the form of feasible ranges.  In this paper we 
calculate exact decycling numbers for the families $C_{3} \times C_{n}$ for $n \geq 3$, 
$C_{4} \times C_{n}$ for $n \geq 4$, 
$C_{n}^{2}$ for $n \geq 4$, and $C_{n}^{3}$ for $n \geq 5$.  Note that $G \times H$ is the 
familiar cartesian product of graphs.  Also, the $n^{th}$ power of $G$, denoted $G^{n}$, is 
defined by $V(G^{n})=V(G)$ and $$E(G^{n})=E(G) \cup \{ uv | u,v \in V(G) 
\enspace \mathrm{ and } \enspace d(u,v) \leq n \}$$

\section{Decycling $C_{m} \times C_{n}$ for $m \leq n$ and $m=3,4$ }

The calculation of the decycling number of the cartesian product of cycles $C_{m} \times C_{n}$ for
$m \leq n$ is mentioned in several papers as an important open problem (for instance [1]). We will
need the following important theorem from [2] both in this section and throughout the paper:

\begin{theorem}
If $G$ is a connected simple graph with maximum degree $\Delta$ then
$$\nabla(G) \geq \frac{|E(G)| - |V(G)| + 1}{\Delta - 1}$$
\end{theorem}

\bigskip

\begin{corollary}
$\nabla(C_{3} \times C_{n}) \geq n+1$.
\end{corollary}

\begin{theorem}
$\nabla(C_{3} \times C_{n}) = n+1$
\end{theorem}

\noindent \textbf{\textit{Proof:}} We need only show that $\nabla(C_{3} \times C_{n}) \leq n+1$. Let
the three canonical $n$-cycles of the graph be $A_{i}$ $1 \leq i \leq 3$ in consecutive cyclic order.
Let the $n$ canonical 3-cycles be $B_{i}$ $1 \leq i \leq n$ in consecutive cyclic order.  Let
$G = C_{3} \times C_{n}$.  Let $C \subseteq V(G)$ be a set of vertices $v_{i}$ such that $v_{i}$ is 
on $B_{i}$ and $v_{i}$ and $v_{i+1}$ are not on the same $A_{j}$ for any $i$. Then $|E(G-C)| = 
6n - 4n = 2n$ and $|V(G-C)| = 3n - n = 2n$.  Finally, it is clear that $G-C$ is connected, so that it is 
therefore unicyclic.  Thus $C \cup \{v\}$, where $v$ is any vertex on the one remaining cycle of $G-C$, is
a decycling set of $G$ with cardinality $n+1$. \hfill $\square$

\bigskip

\begin{theorem}
$\nabla(C_{4} \times C_{n})=\left \lceil \frac{3}{2}n \right \rceil$ for $n \geq 4$
\end{theorem}

\bigskip

\noindent \textbf{\textit{Proof:}} Let $S$ be a decycling set for $C_{4} \times C_{n}$.  Let the 
$n$ canonical cubes which constitute the graph be $q_{i}$ for $1 \leq i \leq n$ in consecutive cyclic 
order.  Let $n_{i}=|V(q_{i}) \cap S|$.  Note that $\nabla(q_{i})=3$ and so
$$3n \leq \sum_{i=1}^{n}n_{i}=2|S|$$
and it follows easily that $|S| \geq \left \lceil \frac{3}{2}n \right \rceil$.

Note that in the following diagrams we draw $C_{4} \times C_{n}$ in its standard toroidal embedding with 
the torus represented as a rectangle with appropriate identifications.  We circle vertices which are members
of a decycling set, and we draw edges of the subgraph induced by the complement of the decycling set more
thickly than the edges which are deleted when the decycling set is deleted.  Figure 1 shows a decycling
set for $C_{4} \times C_{4}$ of cardinality 6, from which it follows that $\nabla (C_{4} \times C_{4})=6$:

\begin{figure}[h]
\begin{center}
\caption{}
\includegraphics{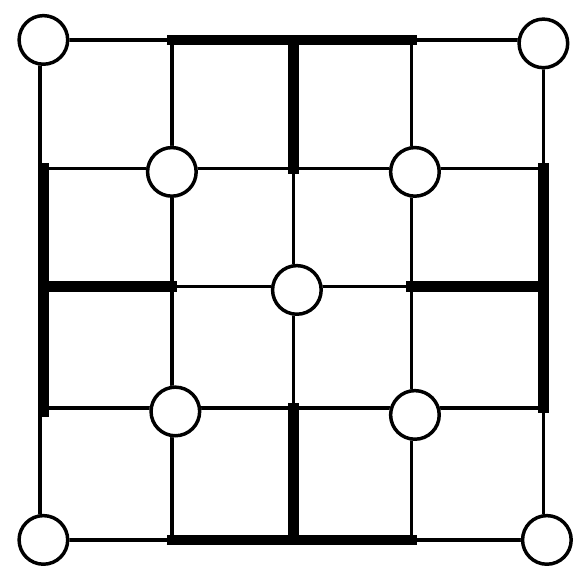}
\end{center}
\end{figure}

\bigskip

\noindent For $C_{4} \times C_{2k}$ we augment this graph by inserting $k-2$ copies of the cylinder
shown in Figure 2.

\begin{figure}[h]
\begin{center}
\caption{}
\includegraphics{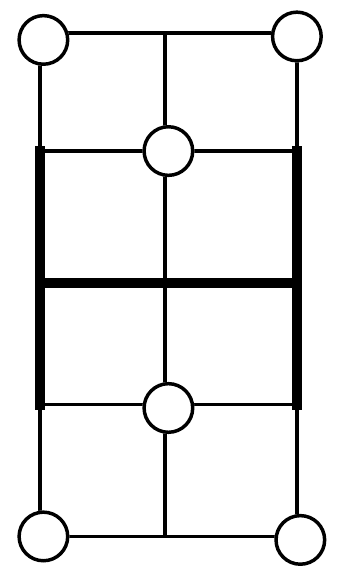}
\end{center}
\end{figure}

The odd case is similar.  For $C_{4} \times C_{5}$ we have the decycling set shown in Figure 3, and then
for $C_{4} \times C_{2k+1}$ we augment by inserting $k-2$ copies of the cylinder shown in Figure 2. 
It is 
easily seen that these decycling sets have the requisite cardinality. \hfill $\square$ 

\begin{figure}[h]
\begin{center}
\caption{}
\includegraphics{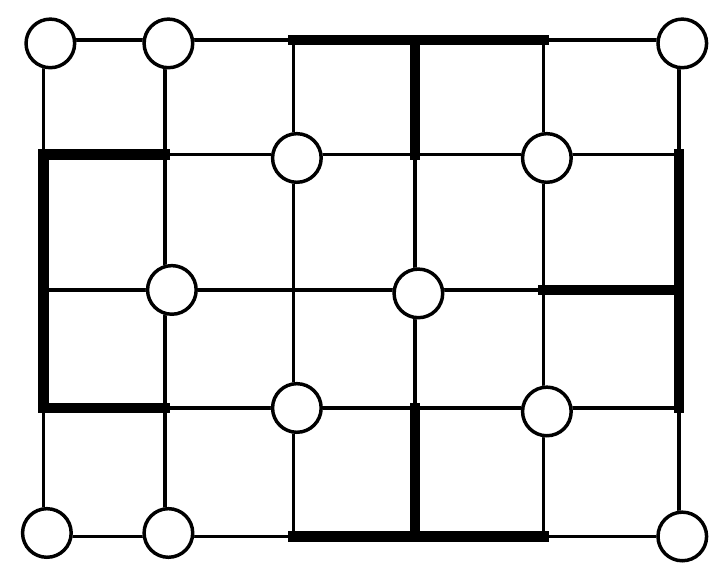}
\end{center}
\end{figure}

\section{Decycling $C_{n}^{2}$ for $n \geq 4$}

The following is a consequence of Theorem 1:

\bigskip

\begin{corollary}
$\nabla \left (C_{n}^{2}\right ) \geq \left \lceil \frac{n+1}{3} \right \rceil$
\end{corollary}

\bigskip

\begin{theorem}
$\nabla \left (C_{n}^{2} \right ) = \left \{ 
\begin{matrix}
&\left \lceil \frac{n+1}{3} \right \rceil & \enspace & n \nequiv 2 \pmod{3}\\

\\

&\left \lceil \frac{n+1}{3}\right \rceil +1 & \enspace & n \equiv 2 \pmod{3}
\end{matrix}
\right. $
\end{theorem}

\bigskip

\noindent \textbf{\textit{Proof:}} Note that for $n \nequiv 2 \pmod{3}$ the lower bound follows from 
Corollary 2.  Let the vertices of $C_{n}^{2}$ be numbered cyclically from 0 to $n-1$.  Note that in the 
remainder of the paper we assume this numbering without special mention.

\bigskip

\noindent \underline{Case 1:} $n \equiv 0 \pmod{3}$

\bigskip

\noindent Let $S=\{0,3, \dots , n-3, n-1 \}$.  The subgraph induced by $V\left(C_{n}^{2}\right)-S$ is 
the path $\{1, 2, 4, 5, \dots , n-5, n-4, n-2 \}$ and hence $S$ decycles $C_{n}^{2}$.  Furthermore
$|S| = \frac{n}{3}+1 = \left \lceil \frac{n+1}{3} \right \rceil$.

\bigskip

\noindent \underline{Case 2:} $n \equiv 1 \pmod{3}$

\bigskip

\noindent Let $S=\{0,3, \dots,n-1 \}$. Then the subgraph induced by $V\left(C_{n}^{2}\right)-S$ is the path
$\{1, 2, 4, 5, \dots, n-3, n-2 \}$ and $|S| = \frac{n-1}{3}+1 = \left \lceil \frac{n+1}{3} \right \rceil$

\bigskip

\noindent \underline{Case 3:} $n \equiv 2 \pmod{3}$

\bigskip

\noindent Let $S=\{0,3,6, \dots , n-2, n-1\}$. It is easy to see that $S$ decycles $C_{n}^{2}$ and 
hence that $\nabla(C_{n}^{2}) \leq \left \lceil \frac{n+1}{3} \right \rceil +1$.  Now suppose we
have a decycling set $S$ with $|S| = \left \lceil \frac{n+1}{3} \right \rceil = \frac{n+1}{3}$.
We first show that $S$ cannot contain two consecutive elements of $V(C_{n}^{2})$.  Suppose by way
of contradiction that $\{0,1\} \subset S$.  Note that $|S-\{0,1\}|=\frac{n-5}{3}$ and that 
$\left | V \left (C_{n}^{2} \right ) - \{0,1\} \right | = n-2$.  Let $B_{i}=\{i, i+1, i+2\}$ and
let $n_{i}=|B_{i} \cap S|$.  Then since each consecutive three elements must contain at least one 
element of $S$ and each element of $S$ is in three of the $B_{i}$'s we have
$$n \leq \sum_{i=0}^{n-1}n_{i} = 3 |S| = n+1$$
\noindent However, $n_{-1}, n_{0} \geq 2$, so there must be some $j$ with $ 2 \leq j \leq n-2$ with 
$n_{j}=0$.  This is a contradiction to the assumption that $S$ decycles $C_{n}^{2}$.

Now, let $i \in V \left (C_{n}^{2} \right )$.  Then $\{i+1, i+2 \} \cap \left (V \left (C_{n}^{2} \right )
-S \right) \neq \O$ and $\{i-1, i-2 \} \cap \left (V \left (C_{n}^{2} \right )
-S \right) \neq \O$.  Hence the degree of $i$ in the subgraph induced by $V \left (C_{n}^{2} \right )
-S$ is at least 2.  Thus that subgraph contains a cycle, and this contradicts the assumption that
$S$ decycles $C_{n}^{2}$. \hfill $\square$

\section{Decycling $C_{n}^{3}$ for $ n \geq 5$}

In this case it turns out that the bound on $\nabla(C_{n}^{3})$ 
given by Theorem 1 is too low in all but a small finite 
number of cases.  The result is:

\bigskip

\begin{theorem}
$$\nabla(C_{n}^{3})= \left \{
\begin{matrix}
&\frac{n+2}{2} & \enspace & n \equiv 0 \pmod 2\\
\\
&\frac{n+1}{2} & \enspace & n \equiv 1 \pmod 4\\
\\
&\frac{n+3}{2} & \enspace & n \equiv 3 \pmod 4\\
\end{matrix} \right.
$$
\end{theorem}

\bigskip

\noindent \textbf{\textit{Proof:}} As above we let $B_{i}=\{i,i+1,i+2,i+3\}$ and $n_{i}=|B_{i} \cap S|$ for 
a decycling set $S$.

\bigskip

\noindent \underline{Case 1}: $n \equiv 1 \pmod 4$

\bigskip

\noindent First of all, $S=\{0,1,2,5,6,9,10,\dots,n-4,n-3\}$ decycles $C_{n}^{3}$ and has the advertised 
cardinality.  Now suppose that $S$ decycles $C_{n}^{3}$ and $|S| = \frac{n-1}{2}$.  Note that $n_{i} \geq 2$, 
for otherwise the isomorph of $K_{4}$ induced by $B_{i}$ is not decycled.  Hence we have 
$$2n \leq \sum_{i=0}^{n-1}n_{i}=4|S|=2(n-1)$$
which is a contradiction.

\bigskip

\noindent \underline{Case 2}: $n \equiv 0 \pmod 2$

\bigskip

Let $S=\{0,1,2,4,6,\dots,n-2$.  Clearly $S$ decycles $C_{n}^{2}$ and has the appropriate cardinality. Now
suppose $S$ decycles $C_{n}^{3}$ and $|S|=\frac{n}{2}$.  As before, $n_{i} \geq 2$ for $0 \leq i \leq n-1$, 
and so
$$2n \leq \sum_{i=0}^{n-1}n_{i}=4|S| = 2n$$
Thus $n_{i}=2$ for $0 \leq i \leq n-1$.  Hence both $(B_{i}-\{i\}) \cap (V(C_{n}^{3})-S)$ and 
$(B_{i-3}-\{i\}) \cap (V(C_{n}^{3})-S)$ are nonempty.  This means that every vertex in the subgraph
induced by $V(C_{n}^{3})$ has degree at least 2, and thus this subgraph contains a cycle, which is a
contradiction.

\bigskip

\noindent \underline{Case 3}: $n \equiv 3 \pmod 4$

\bigskip

Let $S=\{0,1,2,4,5,8,9,\dots,n-4,n-3\}$.  This has the requisite cardinality and decycles the graph.
Now suppose $S$ decycles $C_{n}^{3}$.  If $|\{i, i+1, i+2\} \cap S| \leq 2$ for $0 \leq i \leq n-1$ then consider
$j \in V(C_{n}^{3})-S$.  Note that $j$ is adjacent to a vertex in each of $\{j+1, j+2, j+3\}-S$ and
$\{j-3, j-2, j-1\}-S$.  This, as before, means that every vertex in the subgraph induced by $V(C_{n}^{3}-S)$
has degree at least 2, and thus that $S$ does not decycle $C_{n}^{3}$.  Thus $S$ contains some three consecutive
vertices.  Each of the remaining $\frac{n-3}{4}$ vertex disjoint isomorphs of $K_{4}$ in $C_{n}^{3}$ must contain at 
least two vertices, so
$$|S| \geq 3 + 2 \left ( \frac{n-3}{4}\right) = \frac{n+3}{2}$$
\hfill $\square$

\bigskip

Note that it is fairly easy to show using the methods of the final case in the proof of the previous theorem
that $$|S| \geq k+\frac{n-k}{m+1}(m-1)$$ 

\bigskip

\noindent where $n \equiv k \pmod m+1$
and $0 \leq k \leq n$.  This does not, however, give the exact lower bound in every case, although it is certainly 
better than the lower bound given by Theorem 1.  It also seems to be a more difficult problem to find 
decycling sets of the appropriate cardinalities.  However, these methods might work, with appropriate modification, 
to calculate $\nabla (C_{n}^{m})$ for $m > 3$.

\section{References}

\begin{enumerate}

\item Bau, S. and Beineke, L.  The decycling number of graphs. Australas. J. Combin. 25(2002) pp. 285-298.

\bigskip

\item Beineke, L. and Vandell, R. Decycling graphs. J. Graph Theory. 25(1996) 59-77.

\bigskip

\item Ellis-Monaghan, J. A., Pike, D.A., and Zou, Y. Decycling of Fibonacci cubes. Australas. J. Combin.
35(2006) 31-40.

\bigskip

\item Karp, R.M. Reducibility among combinatorial problems.  In Complexity of computer computations
(Proc. Sympos. IBM Thomas J. Watson Research Center. Yorktown Heights, NY) pp. 85-103. Plenum, New York, 1972.

\end{enumerate}

\end{document}